\documentclass[final,12pt]{jotart2}

\usepackage{amsthm}

\usepackage{amsmath}
\usepackage{amssymb}
\usepackage{verbatim}

\theoremstyle{proclaim}
\newtheorem{theorem}{Theorem}[section]
\newtheorem{lemma}[theorem]{Lemma}
\newtheorem{corollary}[theorem]{Corollary}
\newtheorem{proposition}[theorem]{Proposition}
\theoremstyle{fancyproclaim}
\newtheorem*{vartheorem}{ }

\theoremstyle{statement}
\newtheorem*{rem}{Remark}
\newtheorem*{remarks}{Remarks}
\newtheorem{definition}[theorem]{Definition}
\newtheorem{definitions}[theorem]{Definitions}

\theoremstyle{fancystatement}

\numberwithin{equation}{section}

\providecommand{\AMS}{$\mathcal{A}$\kern-.1667em%
\lower.25em\hbox{$\mathcal{M}$}\kern-.125em$\mathcal{S}$}
\begin{document}
\commby{Editor}

\newcommand{\nc}{\newcommand}
\newcommand{\rnc}{\renewcommand}

\renewcommand{\Re}{\operatorname{Re}}
\renewcommand{\Im}{\operatorname{Im}}
\renewcommand{\le}{\leqslant}
\renewcommand{\ge}{\geqslant}
\renewcommand{\epsilon}{\varepsilon}
\newcommand{\Ker}{\operatorname{Ker}}
\newcommand{\inter}{\operatorname{int}}
\newcommand{\spect}{\operatorname{spec}}
\newcommand{\cA}{\mathcal A}
\newcommand{\cF}{\mathcal F}
\newcommand{\cH}{\mathcal H}
\newcommand{\cD}{\mathcal D}
\newcommand{\CF}{\mathcal F}
\newcommand{\CX}{\mathcal X}
\newcommand{\cW}{\mathcal W}
\newcommand{\tr}{\mathbf T}
\newcommand{\La}{\Lambda}
\newcommand{\wt}{\widetilde}
\newcommand{\wh}{\widehat}
\newcommand{\wV}{\widehat V}
\newcommand{\Om}{\Omega}
\newcommand{\Omin}{\Om_{\text{int}}}
\newcommand{\Omext}{\Om_{\text{ext}}}
\newcommand{\lain}{\lambda_{\text{int}}}
\newcommand{\laext}{\lambda_{\text{ext}}}
\newcommand{\al}{\alpha}
\newcommand{\om}{\omega}
\newcommand{\tht}{\theta}
\newcommand{\si}{\sigma}
\newcommand{\bomin}{\overline{\Om}_{\text{int}}}
\newcommand{\bomout}{\overline{\Om}_{\text{ext}}}
\newcommand{\bomext}{\bomout}
\newcommand{\X}{\mathcal X}
\newcommand{\prt}{\partial}
\newcommand{\De}{\Delta}
\newcommand{\BC}{\Bbb C}
\newcommand{\BD}{\Bbb D}
\newcommand{\BT}{\Bbb T}
\newcommand{\BR}{{\Bbb R}}
\newcommand{\A}{\mathbf{A}}
\newcommand{\B}{\mathbf{B}}
\newcommand{\C}{\mathbf{C}}

\newcommand{\cC}{{\mathcal C}}
\newcommand{\CB}{\mathcal B}
\newcommand{\BZ}{\Bbb Z}
\newcommand{\CE}{\mathcal E}
\newcommand{\CH}{\mathcal H}
\newcommand{\cL}{\mathcal L}
\newcommand{\CS}{\mathcal S}
\newcommand{\CN}{\mathcal N}
\newcommand{\CO}{\mathcal O}
\newcommand{\CW}{\mathcal W}
\newcommand{\CV}{\mathcal V}
\newcommand{\cU}{\mathcal U}
\newcommand{\eps}{\epsilon}
\newcommand{\ga}{\gamma}
\newcommand{\n}{\|}

\newcommand{\Ga}{{\Gamma}}
\newcommand{\de}{\delta}
\newcommand{\la}{\lambda}
\newcommand{\be}{\beta}
\newcommand{\ka}{\varkappa}
\newcommand{\hde}{{\CH(\delta)}}
\newcommand{\hdest}{{\CH(\delta^T)}}
\newcommand{\oGa}{\overline{\Gamma}}
\newcommand{\hdt}{{\CH(\delta^\tr)}}
\newcommand{\Ein}{E^2(\Omin)}
\newcommand{\EinX}{E^2(\Omin,X)}
\newcommand{\EinU}{E^2(\Omin,U)}
\newcommand{\BEinU}{E^2(\bomin,U)}
\newcommand{\EinY}{E^2(\Omin,Y)}
\newcommand{\EomX}{E^2(\Om,X)}
\nc{\whMz}{{\wh M_z}}
\newcommand{\EoutY}{E^2(\Omext,Y)}
\newcommand{\EoutU}{E^2(\Omext,U)}
\newcommand{\EextY}{E^2(\Omext,Y)}
\newcommand{\EextU}{E^2(\Omext,U)}
\newcommand{\BEextY}{E^2(\bomext,Y)}
\newcommand{\BEextU}{E^2(\bomext,U)}
\newcommand{\BEoutU}{E^2(\bomout,U)}
\newcommand{\eout}{E^2(\Omext)}
\newcommand{\ein}{E^2(\Omin)}
\newcommand{\Hol}{\operatorname{Hol}}
\newcommand{\Range}{\operatorname{Range}}
\renewcommand{\phi}{\varphi}
\newcommand{\spec}{\operatorname{spec}}
\newcommand{\dist}{\operatorname{dist}}
\newcommand{\clos}{\operatorname{clos}}
\newcommand{\const}{\operatorname{const}}
\newcommand{\Mztr}{M_z^\tr}
\newcommand{\UAC}{U_{A,J}}
\newcommand{\UABST}{U_{A^*,B_*}}
\newcommand{\Mztrout}{M_{z,out}^\tr}
\newcommand{\Mztrbout}{M_{z,\overline{\text{ext}}}^\tr}
\newcommand{\appr}{\overset\smile\to\smallfrown}
\newcommand{\rstr}{\lceil}
\newcommand{\sm}{\setminus}
\newcommand{\lga}{L^2_\ga(\BR_-,R)}
\newcommand{\pt}{\partial}
\newcommand{\Ao}{A_{0}}
\newcommand{\Aoo}{A_{00}}
\newcommand{\Aka}{A_\ka}
\nc{\WoAB}{\overset{o}{W}_{A,B}} \nc{\WAB}{W_{A,B}}
\nc{\WoABSTR}{\overset{o}{W}_{\A,\,\B}} \nc{\WABSTR}{W_{\A,\,\B}}
\nc{\WhABSTR}{\wh W_{\A,\,\B}}
\newcommand{\defn}{\overset\text{\rm def}\to=}
\newcommand{\Hnty}{H^\infty}
\nc{\wtcX} {{\wt\CX}} \nc{\ess}{ \text{ess}} \nc\defin {\overset
{\text {\rm def} }{=}}

\nc\beqn{\begin{equation}} \nc\neqn{\end{equation}}
\nc{\beqnay}{\begin{eqnarray}}   \nc{\neqnay}{\end{eqnarray}}
\nc{\beqnays}{\begin{eqnarray*}} \nc{\neqnays}{\end{eqnarray*}}
\nc{\barr}{\begin{array}}        \nc{\narr}{\end{array}}
\nc\nn{\nonumber}


\title[
Nagy--Foia\c{s} type functional models in parabolic domains
]{
Nagy--Foia\c{s} type functional models
of nondissipative operators
in parabolic domains
}
\author[Dmitry V. Yakubovich]{Dmitry V. Yakubovich$\;{}^{1)}$}
\addtocounter{footnote}{1}
\footnotetext{This work was finished under the support of the
Ram\'on and Cajal Programme (2002) by the Ministry of the Science
and Technology of Spain, the FEDER and the MEC grants MTM2004-03822,
MTM 2005 - 08350-C03-01}

\address{\noindent Dmitry V. Yakubovich, \newline
\phantom{DDe} Dept. of Mathematics,
Universidad Aut\'onoma de Madrid, Spain}
\email{dmitry.yakubovich@uam.es}
\begin{abstract}
A functional model for nondissipative unbounded perturbations
of an unbounded self-adjoint operator on a Hilbert space
$X$ is constructed. This model is analogous to the Nagy--Foias model
of dissipative operators, but it is linearly similar and not unitarily
equivalent to the operator. It is attached to a domain of parabolic type,
instead of a half-plane.
The transformation map from $X$ to the model space and the analogue of the
characteristic function are given explicitly.

All usual consequences of the Nagy--Foias construction (the $H^\infty$ calculus,
the commutant lifting, etc.) hold true in our context.
\end{abstract}

\maketitle

\section*{INTRODUCTION}


The paper is devoted to the construction of a functional model
of
non-dissipative linear operators of the form
\beqn
\label{intro def A}
A=A_0+i\psi(A_0)F\psi(A_0),
\neqn
where $A_0$ is a self-adjoint unbounded linear operator
on a Hilbert space $X$, $F$ is bounded on $X$
(not necessarily self-adjoint) and
$1\le \psi(x)\le K(\sqrt{|x|}+1)$ for $x\in \cD(\psi)$.
We assume that either
$\psi$ is defined on $\BR$
and is \textit{even} or it is defined on $[0,+\infty)$.
The spectrum of $A$ is contained in an unbounded parabolic domain, which is symmetric
with respect to the real axis.
A precise definition of the unbounded operator $A$ and
precise conditions on $\psi$ will be given later. We remark that
this definition is a particular case of that of
\cite{Mark_Mats}.

All Hilbert spaces in this paper are assumed to be
complex and separable. If $X_1, X_2$ are Banach spaces,
$\CB(X_1,X_2)$ will denote the set of all bounded operators from
$X_1$ to $X_2$.

It is widely recognized that
for understanding the spectral structure of an operator,
the method of functional models is one of
the most useful tools. The original Nagy--Foia\c{s} functional model
exists for contractions and dissipative operators, and is attached,
correspondingly, to the unit disc or to the upper half-plane.

This work can be considered as a continuation of
author's paper \cite{YmodAandA}, where
a linearly similar variant of the Sz.-Nagy---Foia\c{s} model was
suggested and studied.
A general scheme for constructing such kind of models was presented and
several concrete examples were given.
This model has many points in common with the original
Sz.-Nagy---Foia\c{s} model, but also has important differences.
In particular, depending on the operator,
it is constructed in a rather general domain in the complex plane
(bounded or unbounded) and not only in a disc or a half-plane.
We will comment on other differences later on.
A related functional model was constructed by A. Tikhonov in \cite{Tikh1};
see also his subsequent works \cite{Tikh2}, \cite{Tikh3}.
Tikhonov's model, in fact, is closer to the model by Naboko
\cite{Naboko}, which, in particular, made it possible to develop a
stationary scattering theory in the non-selfadjoint context.
We remark that, in general, the function theory that
appears in the Nagy--Foia\c{s} model is studied much better than
analytic questions that arise from the model by Naboko.

Operators of the form \eqref{intro def A} frequently appear in applications.
Namely, suppose that
$A_0$ is a selfadjoint elliptic operator
with regular boundary conditions in $L^2(G)$, where
$G$ is a  bounded domain in $\BR^n$
with smooth boundary.
Take $\psi(x)=1+|x|^\al$, where $0<\al\le \frac 12$. Then
$L\defin A-\Ao$ has the desired form $L=i\psi(\Ao)F\psi(\Ao)$, with
a bounded $F$, iff $(I+|\Ao|)^{-\al}L$ is bounded from
$\cD\big((I+|\Ao|)^{\al}\big)$ to $L^2(G)$.
Note that $\cD\big((I+|\Ao|)^{\al}\big)$ is a kind of Sobolev class in $G$.
Typically, $L$ can be a differential operator of
order less or equal than $4\al m$, where
$2m$ is the order of $\Ao$.
The same is true for elliptic operators on closed manifolds.
We refer to \cite{Agranovich}, \cite{Triebel} and others
for details.

The completeness of eigenvectors of operators
of a related class was established by Keldysh, see
\cite{Keldysh}.

As is known, the spectrum of $A$ lies in a suitable
parabolic domain, see
\cite{Markus} and others; the boundary of this domain is
called sometimes the ``Carleman parabola''.
In this work, we apply the
general scheme of \cite{YmodAandA}
and construct a Sz.-Nagy---Foia\c{s} type functional model
of operator $A$ in a parabolic domain of this type.

In fact, two closely related models
were considered in \cite{YmodAandA}, and we will write down both
these models of $A$
explicitly. They were called in
\cite{YmodAandA} \emph{the quotient model} and
\emph{the resolvent model}.
 It turns out that these constructions
have some points in common with the control theory, in particular,
with the theory of $L^2$ well-posed systems, which was developed in works
by Salamon, Curtain, G. Weiss, Staffans and others.
In particular, our models are not unique, and their choice
depends on the inclusion of our operator in a triple $(A,B,C)$, which is an
abstract analogue of linear control system.
Here we adopt the systems theory terminology
and slightly change the terminology of \cite{YmodAandA}.
We will call here the quotient model
\textit{the control model} and the resolvent model
\textit{the observation model}.
In \S1 and \S5,  these terms will be explained.

Let us describe briefly the control model
(whose connection with
the original setting by Nagy and Foia\c{s}
is more transparent).
Recall first the definition of
\emph{the Smirnov class} $\Ein$. It
consists of all functions $f$ analytic in $\Omin$
such that
$\sup_n\int_{\pt \Om_n}|f|^2\,|dz|<\infty$
for a sequence of domains
$\Om_1\subset \Om_2\subset\dots\subset \Om_n\subset\dots$ with
rectifiable boundary and with $\cup_n\Om_n=\Omin$. We refer to
\cite{Dur} for the properties of the Smirnov classes
$E^p(\Omin)$ and their relationship with Hardy classes
$H^p(\Omin)$. The functions in
$\Ein $ have nontangential boundary values a.e. on
$\Ga$. Equipped with the norm
$$
\|f\|_{\Ein}\defin \frac 1{2\pi}\int_\Ga |f(z)|^2\,|dz|,
$$
the class $\Ein$ is a Hilbert space.

For an auxiliary Hilbert space $U$,
the elements of the Hilbert functional space $\EinU\defin \Ein\otimes U$
are $U$-valued functions analytic in $\Omin$.
These functions also have
nontangential boundary values almost everywhere \cite{SzNF}.
The norm in $\EinU$ is given by
$$
\|f\|_{\EinU}\defin \frac 1{2\pi}\int_\Ga \|f(z)\|_U^2\,|dz|.
$$
The space $\EinU$ can be interpreted either as a
closed subspace of $L^2(|dz|,U)$ or as a space of
$U$-valued analytic functions in $\Omin$. We will use both
interpretations.

We need some more definitions.
Let $U$, $Y$ be two auxiliary Hilbert spaces and $\Omin$ a domain in $\BC$
with piecewise smooth boundary $\Ga$.
Let $\de$ be a function in $\Hnty(\Omin,\CB(Y,U))$.
We call $\de$ \emph{admissible in} $\Omin$ if
there is a constant $\eps>0$ such that
$\|\de(z) y\|\ge \eps \|y\|$, $y\in Y$ for
a.e. $z\in \Ga$.
This function is \emph{two-sided admissible in} $\Omin$
if $\de(z)^{-1}$ exists for a.e.
$z\in \Ga$ and $\|\de^{-1}\|\le C$ a.e. on $\Ga$.
Note that the functions in
$\Hnty(\Omin,\CB(Y,U))$ have
nontangential limits in the strong operator topology
a.e. on $\Ga$
(see \cite[\S V.2]{SzNF}).
If $\de$ is admissible, then
the space $\de \EinY$ is a closed subspace of $\EinU$.

For a holomorphic function $f$ in some domain in $\BC$, we set
$$
M_zf(z)\defin zf(z),
$$
so that $M_z$ is the operator of multiplication
by the independent variable. In general, for any function $\eta$,
we denote by $M_\eta$ the
multiplication operator by $\eta$.

Put
$$
\phi(x)=\psi^2(x).
$$
For $0<\mu<\infty$ and $0<R<\infty$, consider  a parabolic-type domain
\beqn 
\Om_\mu^{\text{int}}=\big\{z=x+iy\in\BC:
x\in\inter\cD(\phi),|y|<\mu\phi(x)
\big\},
\neqn 
We set
\beqn
\label{Om mu R}
\Om_{\mu,R}^{\text{int}}
=\Om_\mu^{\text{int}}\cup B_R(0), \qquad
\Om_{\mu,R}^{\text{ext}}=
\BC\sm\clos \Om_{\mu,R}^{\text{int}},
\neqn
where $B_R(\la)$ stands for the open disc in
the complex plane
of radius $R$ centered at $\la$.
It is known that for certain $\mu$ and $R$,
$\Om_{\mu,R}^{\text{int}}$ contains the spectrum of $A$.
The control model of $A$ is given by
Theorem \ref{last thm}. It asserts
that
\textit{for suitable $\mu$ and $R$,
the operator $A$ is similar to the (unbounded) operator
of multiplication by the independent variable
on the quotient space
$$
E^2(\Omin,X)/\de\cdot E^2(\Omin,X),
$$
where
$\Omin=\Om_{\mu,R}^{\text{int}}$ and
$\de$ is a two-sided admissible $\Hnty$ function on $\Omin$ with values in
$\CB(X)$.
}
Function $\de$ plays the role of
the Nagy--Foia\c{s} characteristic function.
It will be given below by an explicit formula.

The model we get is, in fact, an analogue of a
$C_{00}$ type model in the domain $\Om_\text{int}$, which has no
absolutely continuous part corresponding to the boundary curve.

We will derive a few corollaries from our results.
In particular, one can assert that there exists an
unbounded normal dilation
of $A$ (up to similarity), whose spectrum
lies on $\pt \Omin$. See Corollary \ref{cor2} of Theorem \ref{last thm}.

Before formulating the control model,
we prove Theorem \ref{thm observ},
which gives the observation model of $A$.
These two models of $A$ are equivalent.
All necessary definitions will be given below.

It is very interesting to compare our result with results
by Putinar and Sundberg \cite{Put_Sund}
and Badea, Crouzeix, Delyon
\cite{Bad_Crou_Del}
(see also \cite{Crouz_Del_sector}, \cite{Crouz_parb}).
The results of \cite{Put_Sund} imply that
for \textit{any bounded operator} $A$ on a Hilbert space
and a convex domain $\Omin$ such that the numerical range of
$A$ is contained in its closure, one
can find a dilation of $A$ to
an operator similar to a normal one, whose spectrum
is contained in the boundary of  $\Omin$.
The same holds true for the case
of an unbounded $A$, if
$\Omin$ is a sector, see
\cite{leMe3}.
If an analogous result were true for a general unbounded operator
and a general convex domain $\Omin$,
it would give  a better domain $\Omin$ than our results
for the case when $\Ao$ is bounded from below and
$\|F\|=\|F\|_\text{ess}$.
On the other side, if the spectrum of $\Ao$
is unbounded from above and from below, the numerical range of
$A$ can be the whole complex plane, and
the approach of these papers
does not apply. We also remark that in these works, no expression for
a characteristic function was given, and that our methods are completely
different.

We can mention also the works \cite{Arlinsky} and others by Arlinsky,
where characteristic functions of sectorial operators have been investigated.

Our approach is based on \emph{the duality} between
observation models of $A$ and $A^*$ with respect to
a two-sided admissible function $\de$.
This notion was introduced and studied in \cite{YmodAandA}.
Once dual observation systems $(A,C)$ and $(A^*,B^*)$ are found,
they give rise immediately to dual observation models of $A$ and $A^*$.
In order to prove this duality with respect to $\de$, one has to
find auxiliary operators $B$ and $C$ such that $\de$ and the transfer function of
system $(A,B,C)$ are related by a certain algebraic identity. In our case, we are
able to give such $B$ and $C$ explicitly.

A serious disadvantage of our results resides in the fact that the
values of the generalized characteristic function $\de$
are infinite-dimensional operators and not matrices.
In our setting it is inevitable
(because the dimensions of eigen-spaces of $A$
need not be uniformly bounded).
It can be shown that in some important cases,  the characteristic
function $\de$ has a scalar multiple. This will be discussed
elsewhere.

It seems that if more conditions on $\Ao$ are imposed, then one can obtain
a finite-dimensional model of the same type.

The plan of the exposition is as follows. In Section 1, the observation models of
$A$ and $A^*$ and a duality result will be formulated.
In Section 2, we prove the boundedness of the similarity transformation
$\CO_{A,C}$, which goes from $X$ to the observation model space.
In Section 3, we give more background on the duality and formulate the abstract
result from \cite{YmodAandA} that will be used.
In Section 4, we finish the proofs of our results on observation models.
In Section 5, the control model of $A$ will be introduced, and it will
be explained how to pass to it from the observation model.
In the end of this Section, one can encounter some corollaries and a discussion.
Finally, in Section 6 
we prove some auxiliary geometric lemmas that have been
used earlier.


\section{AN OBSERVATION MODEL AND A DUALITY RESULT}

\subsection{Abstract observation systems and almost diagonalizing transform}
We will have to reproduce some notions and results from \cite{YmodAandA},
which will be used here.

In what follows, we will consider linear systems $(\A,\B,\C)$ of
possibly unbounded operators, where $\A$ acts on
a space $X$, $\B$ acts from \textit{an input space}
$U$ to $X$ (or to a larger Hilbert space)
and $\C$ goes from $X$ (or its dense linear subset) to
\textit{an output space} $Y$.
Spaces $X$, $U$, $Y$ are Hilbert.
A pair $(\A,\C)$ of operators as above will be called
\textit{an observation system} and
a pair $(\A,\B)$  will be called
\textit{a control system}.
Despite the parallelism with the infinite dimensional linear
systems theory, in this abstract setting
we do not need the assumption that
the set $\Re\si(A)$ is bounded from above.

We will use bold letters when discussing the
abstract constructions of functional models and
usual ones when referring to our concrete operator
$A$, given by \eqref{intro def A}, and corresponding
auxiliary operators $B$ and $C$, which will be defined below.

\begin{definition}   
A pair of operators  $(\A,\C)$ (possibly, unbounded) will be called
\textit{an abstract observation system} if

\nr{1} $\A$ is a closed densely defined operator on a Hilbert space $X$
with nonempty field of regularity
$\rho(\A)=\BC\sm\si(\A)$;

\nr{2} $\C:\cD(\C)\to Y$, where $\cD(\C)=\cD(\A)\subset X$ and
$\C$ is bounded in the graph norm
$\|x\|_G\defin \big(\|x\|^2+\|\A x\|^2\big)^{1/2}$ in $\cD(\A)$.
\end{definition}

With every abstract observation system
$(\A,\C)$ we associate the transform
$\CO_{\A,\,\C}$, defined by
$$
\CO_{\A,\,\C}x(z)=\C(zI-\A)^{-1}x, \qquad x\in X,\, z\in \rho(\A).
$$
This map acts from $X$ to the space of
$Y$-valued functions analytic on $\rho(\A)$.

Now let $\Omin$ and $\Omext$ be a pair of open subsets in $\BC$ that have a
common boundary $\Ga$.
In this abstract setting, our requirements are:
\textit{(i)} $\Omin\cap\Omext=\emptyset$; $\BC=\Omin\cup\Omext\cup\Ga$;
\textit{(ii)} $\Ga$ is a finite union of piecewise smooth contours, each of them
homeomorphic to the unit circle or a real line. In the latter case, both ends of the
contour have to go to infinity;
\textit{(iii)} $1/(|z|+1)\in L^2(\Ga,|dz|)$.
If these conditions hold, we will call
the open set $\Omin$ \textit{admissible}.
\begin{definition}
We call an abstract observation system $(\A,\C)$
\textit{admissible with respect to} $\Omin$
if $\si(A)\subset \clos\Omin$ and
operator $\CO_{\A,\,\C}$ is bounded from $X$ to $\EoutY$.

We call an abstract observation system $(\A,\C)$
\textit{exact with respect to} $\Omin$
if it
admissible with respect to $\Omin$ and, moreover,
there is a two-sided estimate
$$
\|\CO_{\A,\,\C}x\|_{\EoutY} \asymp \|x\|, \qquad x\in X.
$$
\end{definition}

The relation with the
theory of well-posed control systems is as follows. Suppose now that
$\A$ is a generator of a bounded $C_0$ semigroup.
Let
$$
\Pi_-=\{z:\Re z<0\}, \qquad
\Pi_+=\{z:\Re z>0\}
$$
be, respectively, the left and the right half-planes.
Consider the
linear continuous time observation system
$$
\dot x(t)=\A x(t), \quad x(0)=x_0, \qquad y(t)=\C x(t).
$$
For any initial value $x_0\in \cD(\A)$, the output $y=y(t)$ is a well-defined
continuous function on $[0,+\infty)$. Denote $y=\wh \CO_{\A,\,\C}x_0$, so that
$\wh \CO_{\A,\,\C}$ is the space--output map. Then
$$
\CO_{\A,\,\C}x_0(z)=\big(\cL \wh \CO_{\A,\,\C}x_0\big)(z), \qquad x_0\in \cD(\A)
$$
for all $z\in \Pi_+$,
where $\cL y(z)=\int_0^\infty e^{-zt}y(t)\,dt$ is the Laplace transform.
It follows that in this case, abstract observation system
$(\A,\C)$ is admissible with respect to
$\Omin:=\Pi_-$ if and only if the inequality
$$
\int_0^\infty \|y(t)\|^2\, dt\le K \|x_0\|^2
$$
holds for some constant $K>0$ and for all initial data $x_0\in \cD(A)$.
System $(A,C)$ is exact with respect to $\Pi_-$ if and only if the two sides
of this inequality are comparable.
There is a close connection between this setting and the definition of a
well-posed output map, see  \cite{Sta-libro}, Theorem 4.4.2.


Now let us return to the general situation of an abstract observation
system and arbitrary admissible domain $\Omin$.

\begin{definition}[\cite{YmodAandA}]
Let $\de\in \Hnty(\Omin,\cL(Y,U))$ be
a two-sided admissible function.
We introduce
\emph{the observation model space}
$\hde$ as the following closed subspace of $\EextY$:
\beqn
\label{def hde}
\hde=\big\{
f\in \EextY: \wt f\defin\de \cdot f|\Ga\in \EinU
\big\}.
\neqn
\end{definition}

We introduce
(possibly, unbounded) operators
$\Mztr$, $j$  on $\hde$ as follows. Put
$$
\cD(j)=\cD(\Mztr)=
\big\{
f\in \hde: \exists c\in Y: zf-c\in \hde
\big\}.
$$
For $f\in\cD(\Mztr)$, the constant $c$ is unique. Therefore,
the operators
\begin{gather*}
j:\cD(\Mztr)\to Y, \qquad \Mztr:\cD(\Mztr)\to \hde, \\
jf\defin c, \qquad \big(\Mztr f\big)(z)\defin zf-c,\; f\in \cD(\Mztr)
\end{gather*}
are well defined. We shall call $\Mztr$ the \textit{operator
of truncated multiplication} on $\hde$.

Let  $\de\in \Hnty\big(\Omin, \CB(Y,U)\big)$ be a two-sided admissible function.
By definition (see \cite{NikBook2}),
 \textit{the spectrum} of $\de$ is the set of points $\la\in \clos\Omin$
such that $\de^{-1}\notin \Hnty\big(\Omin\cap\cW, \CB(Y,U)\big)$ for any
neighbourhood $\cW$ of $\la$. It will be denoted by $\spect \de$.
It is a closed subset of $\clos\Omin$.
Its intersection with $\Omin$ consists of points $\la$ in $\Omin$ such that
$\de(\la)$ is not invertible.


Any function $f\in\hde$ can be viewed as an
analytic function
on $\Omext\cup\Omin\sm\spect\de$.
On $\Omin\sm\spect\de$, we define it by
means of the formula
\beqn
\label{ext f z}
f(z) \defin \de(z)^{-1}\wt f(z),
\neqn
where  $\wt f(z)$ is determined from \eqref{def hde}.

For completeness, we give here the formula for the
resolvent of $\Mztr$.

\begin{proposition}[\cite{YmodAandA}, Props. 1.1 and 2.3]

\nr{i} The operator $\big(\Mztr-\la I\big)^{-1}$
exists and is bounded if and only if
$\la\in \BC\sm \spec\de$.

\nr{ii} Each function
$f\in \hde$ extends analytically to $\BC\sm \spect\de$ according
to the rule $f(\la)\defin j(\la I-\Mztr)^{-1}f$.
For $\la\in\Omin\sm\spect\de$,
this extension coincides with
that defined in \eqref{ext f z}.

\nr{iii} For $\la$ in $\BC\sm\spect\de$,
$$
\big(\Mztr-\la I\big)^{-1}
f(z)=\frac{f(z)-f(\la)}{z-\la}, \qquad f\in\hde.
$$
\end{proposition}

Now this scheme will be concretized in order to give
a precise observation model of
operator \eqref{intro def A}.

\subsection{Conditions on the perturbation}
Remind that the essential norm of $F$
is defined as
$$
\|F\|_\ess\defin\inf \{\|F+R\|: R\in\CB(X) \text{ is compact}\}.
$$
Assume that $\psi:\cD(\psi)\to\BR$ and $A$
satisfy the following.

\

\nr{1} Either $\psi$ is defined on $\BR$
and is \textit{even} or it is defined
on $[0,+\infty)$;

\nr{2} \enspace $\psi$ is a continuous function; moreover,
\enspace $\psi$ \enspace  is of
class $C^1$ on $\cD(\psi)\sm\{0\}$;

\nr{3}
$\psi\ge 1$ on $\cD(\psi)$ and $\psi(x)\to+\infty$ as
$x\to+\infty$;

\nr{4}   $\psi^2$ is concave on $[0,\infty)$;

\nr{5} One has
\beqn
\label{cond 5 psi}
\|F\|_\ess\cdot k_0(\psi)<1, \qquad \text{where }\quad
k_0(\psi)\defin \lim_{t\to+\infty}\frac{\psi^2(t)}{{t}}
\neqn
(it follows from (4) that this limit exists).

\nr{6} If $\cD(\psi)=[0,\infty)$, then $\si(\Ao)\subset [\eps_0,\infty)$ for some
$\eps_0>0$.

\

Notice that condition (5) is
automatically fulfilled whenever either
$F$ is compact or $k_0(\psi)=0$.

We put
$$
\phi(t)\defin\psi^2(t).
$$

\subsection{Precise definition of $A$}
Put $\Aoo=I+|A_0|$; then $\cD(A_0)=\cD(A_{00})$.
We rewrite $A$ in the form
\beqn
\label{def A}
A=\Aoo\big[\Aoo^{-1}\Ao+i(\Aoo^{-1}\psi(\Ao))F\psi(\Ao)\big]
\neqn
and take it for the precise definition of $A$. We set
\beqn
\label{defin A}
\begin{aligned}
\cD(A)&\defin \big\{ x\in \cD(\psi(\Ao)): \\ & \qquad\qquad
\big[\Aoo^{-1}\Ao+i(\Aoo^{-1}\psi(\Ao))F\psi(\Ao)\big]x\in
\cD(\Ao) \big\}
\end{aligned}
\neqn
Notice that operators $\Aoo^{-1}\Ao$ and
$\Aoo^{-1}\psi(\Ao)$ are bounded.

Consider the control system $(A,B,C)$, where
$$
C=i\psi(\Ao), \quad B=\psi(\Ao).
$$
We put $Y=U=X$. Notice that formally, $A=\Ao+L$, where the perturbation $L$ factorizes as
\beqn
\label{L eq BFC}
L=BFC.
\neqn
According to \eqref{defin A}, the pair
$(A,C)$ is an abstract observation system.

\begin{definition}
Let $\A$ be a closed densely defined operator on $X$ with
$\si(\A)\ne \BC$. Take any point $\la\in\rho(\A)$.
We define
\emph{the Hilbert space $X_\la(\A)$} as the
vector space of
formal expressions $(\A-\la I)x$,
where  $x$ ranges over the whole space $X$.
Introduce a Hilbert norm on
$X_\la(\A)$ by setting $\|(\A-\la I)x\|_{X_\la(\A)}\defin \|x\|_X$ for all
$x\in X$. For $x\in \cD(\A)\subset X$, we identify the
element $(\A-\la I)x$ of $X(\A)$ with the element of $X$, given by
the same expression.
\end{definition}

It is clear that by this construction, $X$ becomes a dense subset of
$X_\la(\A)$.
This construction does not depend on the choice of $\la$ in the sense that
for different $\la$'s in $\rho(\A)$,
the corresponding spaces $X_\la(\A)$ are
naturally isomorphic
(and have equivalent norms).
If the exact form of the norm in $X_\la(\A)$ is not important,
then we write $X(\A)$ instead of $X_\la(\A)$.
Observe that $\A$ is a bounded operator from $X$ to $X(\A)$.

\subsection{Observation model of $A$}
Put
\beqn
\label{def mu0}
\mu_0=\frac{\|F\|_\ess}{\sqrt{1-\|F\|_\ess^2k_0(\psi)^2}}.
\neqn
For $\ka\in\BR$, we consider the normal (possibly unbounded) operator
$$
A_\ka=\Ao+i\ka\phi(\Ao).
$$
Now we can formulate the ``observation form'' of our functional model.

\begin{theorem}[An observation model of $A$]
\label{thm observ}
Take any $\mu>\mu_0$.
For $\ka\in \BR$, define
$$
\hfill \de_\ka(z)=
\big[\Aka-zI+i\phi(\Ao)F\big]^{-1}
\big[\Ao-zI+i\phi(\Ao)F\big].  \hfill
$$
Then there exist  $R>0$ and $\ka\in \BR$
such that for the corresponding
function $\de\defin \de_\ka$ and for the domains
\beqn
\label{Omin Omext}
\Omin\defin \Om_{\mu,R}^\text{int},
\qquad
\Omext\defin \Om_{\mu,R}^\text{ext}
\neqn
(see \eqref{Om mu R}) the following statements hold.

\nr{1} $A$ is a closed densely defined operator, and
$\si(A)\subset \Omin$. The pair $(A,C)$ is an exact observation system.

\nr{2} Function $\de$ is in $\Hnty(\Omin,\CB(X))$ and is two-sided admissible;

\nr{3} Operator
$$
\CO_{A,C}:X\to \hde
$$
is an isomorphism that transforms
operator $A$ into the truncated multiplication operator
$\Mztr$ on the observation model space $\hde$. This means that
for any $x\in \cD(A)$,
$\CO_{A,C}\,x$ is in $\cD(\Mztr)$,
\beqn
\label{split prop}
\CO_{A,C}\,Ax=\Mztr\CO_{A,C}\,x
\neqn
and that, moreover,
\beqn
\label{split prop2}
\CO_{A,C}\cD(A)=\cD(\Mztr).
\neqn
\end{theorem}

In fact, we will show that there is $\ka_0>0$ such that
one can take any $\ka$, $|\ka|>\ka_0$ in the above theorem.
The value of $\ka_0$ is given below in \eqref{def ka0}.

The above definition of $\de$ can be understood as follows. It is
clear that $\Ao-zI+i\phi(\Ao)F$ is a bounded operator from $X$ to
$X(\Ao)$. It will be proved in Lemma \ref{lem about de}
that for $z\in\Omin$ (and for $\ka>\ka_0$),
the operator $\Aka-zI+i\phi(\Ao)F$ from $X$ to $X(\Ao)$
is invertible. Hence $\de_\ka$ is a
well-defined bounded operator on $X$ for $z\in \Omin$.

The splitting property \eqref{split prop} is in fact a matter of algebra and
holds true in much more general context, see \cite[Proposition 1.2]{YmodAandA}.

This theorem models the operator $A$ by the operator $\Mztr$ on the model
functional space $\hde$.
As it will be seen in \S\ref{sect contr mod},
this model is  very closely related to the Nagy and Foias model.

\subsection{The duality result}

Let $\de \in \Hnty\big(\Omin, \CB(Y,U)\big)$ be a two-sided admissible
operator function in
an admissible domain $\Omin$.
We orient the
curves that constitute
$\Ga=\prt\Omin$ in such a way that, under the movement
along them, the domain $\Omin$ remain on the left.
Put
\begin{gather*}
\bomin\defin\{\bar z:z\in \Omin\},  \\
\de^\tr(z)=\de(\bar z)^*, \qquad z\in \bomin.
\end{gather*}
Then $\de^\tr$ is a two-sided admissible function in $\bomin$.
We will need the model space
\beqn
\label{def hdest}
\nn
\hdest\defin\big\{
f\in \BEextY:\quad  \de^\tr \cdot f|\pt\bomin\in \BEinU
\big\},
\neqn
which is associated to the function
$\de^\tr$ and the domain $\bomin$.

We start with the following fact.

\begin{proposition}[\cite{YmodAandA}, Prop. 4.2]
For any two-sided admissible function
$\de \in \Hnty\big(\Omin, \CB(Y,U)\big)$,
model spaces $\hde$ and $\hdest$ are dual to each other with respect to the
Hermitian pairing
$$
\langle
f, g
\rangle_\de
\defin
\frac 1{2\pi i}\int_\Ga
\langle
\de(z)f(z), g(\bar z)
\rangle
\; dz.
$$
\end{proposition}

In fact, in our case, $\Omin$ is symmetric with respect to the real line, that is
$\Omin=\bomin$.

\begin{definition}
Suppose we are given a triple of
(possibly unbounded) operators
$(\A,\B,\C)$
and Hilbert spaces $X$, $U$, $Y$, which have the meaning of
the state space, the input space and the output space, respectively.
We say that the triple
$(\A,\B,\C)$ is \emph{a full abstract system}
if the following condition hold.

\nr{1} $(\A,\C)$ is an abstract observation system, whose output space is
$Y$;

\nr{2} $\B: U\to X(\A)$ is a bounded operator.
\end{definition}

Each linear continuous functional in
$X(\A)^*$
can be considered in the same time as a
linear continuous functional on $X$, that is, an element of
$X^*$. In this sense,
$X(\A)^*$ coincides with $\cD(\A^*)$, see, for example,
\cite{Sta-libro}. It follows that
whenever $(\A,\B,\C)$ is a full abstract system,
$(\A^*,\C^*,\B^*)$ is also a full abstract system, with input and output space
interchanged.

\begin{definition}
Suppose $\Omin$ is fixed and $(\A,\B,\C)$ is a full abstract system.
Let $\de \in \Hnty\big(\Omin, \CB(Y,U)\big)$ be two-sided admissible.
We say that \textit{observation systems $(\A,\C)$ and
$(\A^*,\B^*)$ are in duality} with respect to $\de$ if

\nr{1} Operators
$$
\CO_{\A,\C}:X\to\hde, \quad \CO_{\A^*,\B^*}:X\to\hdest
$$
are isomorphisms;

\nr{2} For all $x_1,x_2$ in $X$,
\beqn
\label{eq dual}
\langle x_1, x_2
\rangle_X
=
\langle \CO_{\A, \C}x_1, \CO_{\A^*, \B^*}x_2
\rangle_\de.
\neqn
\end{definition}

In addition to Theorem \ref{thm observ}, we will prove the following result.


\begin{theorem}
\label{thm dual}
Take any $\mu>\mu_0$.
Then there exist  $R>0$ and $\ka\in \BR$
such that all statements of Theorem \ref{thm observ} hold and, moreover,
systems
$(A,-\ka C)$ and
$(A^*,B^*)$ are in duality with respect to the function $\de=\de_\ka$.
It is assumed here that $\Omin$ is defined by \eqref{Omin Omext}.
\end{theorem}

This result implies that the transform $\CO_{A^*, B^*}$
is an isomorphism that converts the action of $A^*$  into the action of
$\Mztr$ on the model space  $\hdest$.

Observation models and control models are closely related.
One of these relations is given below in Lemma \ref{rel ctr obs mods}.
Another one is \cite[Prop. 4.1]{YmodAandA}.

\begin{remarks}
\nr{1} If $\de$ is a two-sided admissible function
and $\CO_{\A,\,\C}$ is an isomorphism of $X$ onto $\hde$,
then we called $\de$ a \emph{generalized characteristic function of
an observation system $(\A,\C)$} in \cite{YmodAandA}. Its determination
is far from unique, to the opposite to
the classical notion of the Nagy---Foia\c{s} characteristic function, which
is essentially unique. In fact, it is easy to see that
$\cH(\de)=\cH(\be\cdot\de)$ for any
function $\be$, which is invertible in the algebra
$\Hnty\big(\Omin, \CB(X)\big)$. Therefore
for any $\be$ of this kind,
$\be\cdot\de$ is a generalized characteristic function of
system $(\A,\C)$ together with $\de$.

We would obtain formally a closer analogue of the
classical  Nagy--Foia\c{s} construction if we required $\de$ to be
two-sided inner. However, it is this freedom of the choice of $\de$
that permits us to give an explicit formula for the generalized
characteristic function of the operator $A$ in study.

In \cite{YmodAandA}, in fact,  we discussed
functional models for a more general class of $*$-admissible functions.

\nr{2} Full systems $(\A,\B,\C)$ such that
observation systems $(\A,\C)$ and $(\A^*,\B^*)$ are in duality are very
special ones. Their consideration is motivated by
our scheme of constructing functional models rather than by
the control theory. If a generalized
characteristic function $\de$ of system
$(\A,\C)$ is fixed, then, by \cite[Prop. 9.3]{YmodAandA},
there is a unique operator $\B$ such that
$(\A,\C)$ and $(\A^*,\B^*)$ are in duality with respect to $\de$.
\end{remarks}

\section{ADMISSIBILITY OF THE OBSERVATION SYSTEM $(A,C)$}


We fix some $\mu>\mu_0$.
From now on, let us also fix a number $r'>\|F\|_\ess$, close to
$\|F\|_\ess$, and a number $k>k_0(\psi)$, close to
$k_0(\psi)$, so that
\beqn
\label{ineq r' k mu}
r'k<1, \qquad \frac {r'}{\sqrt{1-{r'}^2k^2}}<\mu
\neqn
(see \eqref{cond 5 psi} and \eqref{def mu0}).
Take any decomposition $F=F'+F''$  such that
$F''$ is a finite rank operator and
\beqn
\label{nFdash}
\| F' \| < r'.
\neqn
It is possible, because any compact operator in $X$ can be
approximated in norm by finite rank operators.

First let us formulate two technical lemmas,
whose proofs will be given in the last Section.

\begin{lemma}
\label{lem un discs}
There exists  $R_0>0$ \enspace
such that for all $t\in\cD(\phi)$,
\beqn
\label{B subset Om}
B(t, r'\phi(t))\subset \Om^{\text{int}}_{\mu,R_0}.
\neqn
\end{lemma}

\begin{lemma}
\label{lem est int}
Let $\Omin=\Om_{\mu,\,R}^{\operatorname{int}}$
for some positive  $R$, and let $\Ga=\pt \Omin$.
Then
there is a positive constant $K$ such that the inequality
$$
|\psi(x)|^2\int_\Ga \frac{|d\la|}{|x-\la|^2}\le K
$$
holds for all $x\in \si(A_0)$.
\end{lemma}


\begin{lemma}
\label{lem Ao C adm}
The system $(A_0,C)$ is admissible
with respect to the domain $\Omin$.
\end{lemma}

\begin{proof}
By the Spectral Theorem, $\Ao$ is unitarily equivalent
to the operator $\wt\Ao f(t)=tf(t)$, acting on a direct integral
$$
\wtcX\defin \int^\oplus \CX(t)\,d\nu(t),
$$
where $\nu $ is a positive Borel measure on $\si(\Ao)$
(a scalar spectral measure of $\Ao$) and
$\{\wtcX(t)\}$ is a $\nu$-measurable family of Hilbert spaces
\cite{BirmSol}.
The same unitary isomorphism converts $C=i\psi(\Ao)$ into
$\wt C=M_{i\psi}$.
We prove our statement by passing to this model of
the pair $(\Ao, C)$.
For $f=f(t)\in \wtcX$,
\begin{multline}
\|(\CO_{\wt\Ao,\wt C}f)\|_{\EextY}^2
=\int_\Ga \int_\BR
\bigg\|
\frac{\psi(t)f(t)}{t-z}
\bigg\|^2\, d\nu(t)\,|dz| \\
=
\int_{\BR}
\, d\nu(t)
\|f(t)\|^2|\psi(t)|^2
\int_{\Ga}
\frac {|dz|}
{|t-z|^2}                \\                       
\le K
\int_{\BR}
\|f(t)\|^2|\, d\nu(t)
=K \|f\|_\wtcX^2.
\nn
\end{multline}
The inequality is due to Lemma \ref{lem est int}.
\end{proof}

\begin{lemma}
\label{lem norm res}
\nr{i} \enspace
$
\displaystyle
\limsup_{z\in\Om_{\mu}^\text{ext},\,z\to\infty}
\| FC(A_0-zI)^{-1}B\|\le \frac{\|F'\|}{r'}<1;
$

\hskip 2.4cm \nr{ii} \enspace
$
\displaystyle
\limsup_{z\in\Om_{\mu}^\text{ext},\,z\to\infty}
\| C(A_0-zI)^{-1}BF\|\le \frac{\|F'\|}{r'}<1.
$
\end{lemma}


\begin{proof}  
By Lemma \ref{lem un discs}, if $|z-t|< r'\phi(t)$ for
some $t\in\cD(\phi)$, then $z\in\Om^{\text{int}}_{\mu,R_0}$.
It follows that for $z\in \clos\Om^{\text{ext}}_{\mu,R_0}$,
\beqn
\label{est phi z-t}
\phi(t)/|z-t|\le 1/r', \qquad \forall t\in \cD(\phi).
\neqn
Hence
\beqn
\label{le 1r'}
\|
\phi(\Ao)(\Ao-zI)^{-1}
\|\le 1/r'
\neqn
for $z\in \clos\Om^{\text{ext}}_{\mu,R_0}$. Moreover, if we put
$\wt\cA(z)=\phi(\Ao)(\Ao-zI)^{-1}$, then it is easy to check that
\beqn
\label{SOT lim}
\wt\cA(z)^*
\rightarrow 0 \quad \text{as}\quad z\to\infty,\; z\in \clos\Om^{\text{ext}}_{\mu,R_0}
\neqn
in the strong operator topology. (One can apply here the Spectral Theorem
in the same way as in the proof of Lemma \ref{lem Ao C adm},
the Lebesgue dominated convergence theorem and \eqref{est phi z-t}.)
We have
$$
\| FC(A_0-zI)^{-1}B\|
\le
\|F'\phi(\Ao)(\Ao-zI)^{-1}\|+
\|F''\phi(\Ao)(\Ao-zI)^{-1}\|.
$$
The relation
\eqref{SOT lim} and the fact that $F''$ has a finite rank imply that
$\|F''(\Ao-zI)^{-1}\phi(\Ao)\|\to 0$ as
$z\to\infty,\, z\in \clos\Om^{\text{out}}_{\mu,R_0}$.
By applying the estimate $\|F'\|< r'$ and \eqref{le 1r'}, we obtain (i).
Assertion (ii) is obtained similarly.
\end{proof}

From now on, we fix some $\eps>0$ and a
radius $R>R_0$ such that
\beqn
\label{n FC resolv B}
\|FC(A_0-zI)^{-1}B\|\le 1-\eps, \qquad
\|C(A_0-zI)^{-1}BF\|\le 1-\eps
\neqn
for all $z\in\Om_{\mu,R}^{\text{ext}}$.
It is possible due to
Lemma \ref{lem norm res}.
According to \eqref{Om mu R}, we put
$\Omin=\Om_{\mu,R}^{\text{int}}$,
$\Omext=\BC\sm\clos\Omin$.

\begin{definition}
Let $\eta$ be a real Borel function on $\cD(\psi)$ such that
$\eta(t)\ne0$ for all $t\in\cD(\psi)$. We define
the Hilbert space
$X_\eta$ as the set of formal expressions $\eta(\Ao)x$, $x\in X$.
Recall that the self-adjoint operator
$\eta(\Ao)$ is bounded iff $\eta$ is essentially
bounded with respect to the spectral measure of $\Ao$.
We introduce a Hilbert norm on
$X_\eta$ by setting $\|\eta(\Ao)x\|_{X_\eta}\defin \|x\|_X$ for all
$x\in X$. For $x\in \cD\big(\eta(\Ao)\big)\subset X$, we identify the
element $\eta(\Ao)x$ of $X_\eta$ with the element of $X$, given by
the same expression. Notice that if $\eta$ is essentially bounded, then
$X_\eta=\cD\big(\eta^{-1}(\Ao)\big)\subset X$.
\end{definition}

This definition is very close to the
definition of spaces $X_\la(\A)$, given earlier.
In fact, if $\eta>\eps>0$, then $X_\eta=X_0(\eta(\Ao))$.

Consider, in particular, the Hilbert space $X_\psi\supset X$.
Since $B$ is an isometric isomorphism of $X$ onto $X_\psi$,
for any $T\in\CB(X_\psi)$,
\beqn
\label{ident norm R}
\|T\|_{\CB(X_\psi)}=\|B^{-1}TB\|_{\CB(X)}.
\neqn
Recalling the notation from
 \eqref{L eq BFC}, by the first inequality in \eqref{n FC resolv B}
we obtain
\beqn
\label{norm L by resolv}
\|L(A_0-zI)^{-1}\|_{\CB(X_\psi)}\le 1-\eps, \qquad
z\in\clos\Omext.
\neqn

\begin{lemma}
\label{lem on siA and H}
\nr{i} $\si(A)\subset\Omin$.

\nr{ii} For $z\in\clos\Omext$, $(A-zI)^{-1}$ is an isomorphism of
$X_\psi$ onto $X_{\frac{\psi(t)}{|t|+1}}$.

\nr{iii} The identity
\beqn
\label{cU eq Hz cU}
\CO_{\Ao,C}x(z)=H(z) \CO_{A,C}x(z), \qquad x\in X, \, z\in \Omext
\neqn
holds, where
\beqn
\label{def H z}
H(z)=I+C(\Ao-zI)^{-1}BF.
\neqn

\nr{iv} One has $H, H^{-1}\in \Hnty(\Omext, \CB(X))$.

\nr{v} Observation system $(A,C)$ is
admissible with respect to the domain $\Omin$.
\end{lemma}

Notice that $\psi^2(t)\le K\big(|t|+1\big)$ implies that
$X_{\frac{\psi(t)}{|t|+1}}\hookrightarrow X_{1/\psi}$.

\begin{proof}[Proof of Lemma \ref{lem on siA and H}]
The definition
\eqref{defin A} of $\cD(A)$ can be rewritten as
$$
\cD(A)\defin
\big\{
x\in X_{1/\psi}: [\Ao+i\psi(\Ao)F\psi(\Ao)]x\in X
\big\},
$$
where $[\Ao+i\psi(\Ao))F\psi(\Ao)]x$ is understood a priori as an element of
$X_{|t|+1}$. Recall that $X_{1/\psi}=
\cD\big(\psi(\Ao)\big)$.
Hence for all $y\in X_{\frac{\psi(t)}{|t|+1}}$,
the equality
$$
(A-zI)y=
\big(
I+L(\Ao-zI)^{-1}
\big)
(\Ao-zI)y
$$
between elements of $X_\psi$ holds for all $z\notin \si(\Ao)$.
By \eqref{norm L by resolv},
$I+L(\Ao-zI)^{-1}$ is invertible in $X_\psi$ for $z\in \clos\Omext$.
Hence for these $z$, $A-zI$ has a bounded inverse in $X$, given by
\beqn
\label{frm for res A}
(A-zI)^{-1}=(\Ao-zI)^{-1}
\big(
I+L(\Ao-zI)^{-1}
\big)^{-1}
\neqn
(notice that the immersion $X\hookrightarrow X_\psi$ is bounded).
This proves (i). Formula \eqref{frm for res A} also gives (ii).

Similarly, we have
$$
(A-zI)y=(\Ao-zI)\big(I+(\Ao-zI)^{-1}BFC\big)y, \qquad y\in \cD(\psi(\Ao)),
$$
which implies that
$$
(\Ao-zI)^{-1}
=
(I+(\Ao-zI)^{-1}BFC)
(A-zI)^{-1},
$$
where $I+(\Ao-zI)^{-1}BFC$ is a  bounded operator  from $\cD(\psi(\Ao))$ to $X$.
By multiplying this equality by $C$ from the left, we obtain (iii).

Assertion  (iv) follows from the second inequality in \eqref{n FC resolv B}.
At last,  \eqref{cU eq Hz cU}, (iv)
and Lemma \ref{lem Ao C adm} imply (v).
\end{proof}

Assertion (ii) of Lemma implies that $A$ is closed and densely defined.

Notice that the inequality \eqref{norm L by resolv} implies that
\beqn
\label{norm Omext}
\|\big(
I+L(\Ao-zI)^{-1}
\big)^{-1}\|_{\CB(X_\psi)}
\le \eps^{-1}<\infty, \qquad
z\in\Omext.
\neqn

\section{OUTLINE OF THE PROOF OF THEOREM \ref{thm dual}}

Let $(\A,\B,\C)$ be a full abstract system
with the state space $X$, input space $U$
and output space $Y$. Let
$\Phi$ be a holomorphic operator-valued
function on $\rho(A)$ with values in $\CB(U,Y)$.

\begin{definition}
We call $\Phi$ \textit{a transfer function }
of system $(\A,\B,\C)$ if the identity
$$
\Phi(z)-\Phi(w)=\C
\big[
(zI-\A)^{-1}-(wI-\A)^{-1}
\big]\B
$$
holds for all $z,w\in \rho(A)$.
\end{definition}

This definition is standard in the theory of well-posed systems, see
\cite{Sta-libro}. The point is that the difference of the resolvents
of $\A$ in points $z$ and $w$ is a bounded map from $X(\A)$ to
$\cD(\A)$, which implies that the right hand part is
always in $\CB(X)$.
The transfer function of a system is determined uniquely up to adding
an arbitrary operator constant.

We need the following definition from \cite{YmodAandA}.

\begin{definition}
We tell that
a function
$\Phi\in \Hnty(\Omext,\CB(U,Y))$
\textit{corresponds to}
a two-sided admissible function $\de\in \Hnty(\Omin,\CB(Y,U))$
if there is a function $\tau\in \Hnty(\Omin,\CB(U,Y))$
such that
the following two conditions hold:

1) $\Phi|\Omext\in \Hnty(\Omext,\CB(U,Y))$;

2) $\Phi_e=(\de^{-1}+\tau)_i$ a.e. on $\Ga$.
\end{definition}
Our main tool in proving
Theorem \ref{thm dual} will be the following result from
\cite{YmodAandA}.

\begin{vartheorem}[Theorem A]
{\normalfont (see \cite{YmodAandA}, Theorem 9.5) }
\label{thm 9-5}
Let $(\A,\B,\C)$ be a full abstract system and $\Phi$ its transfer function.
Suppose that $\Omin$ is an admissible domain,
$\si(\A)\subset\clos\Omin$, and let $\de$ be a two-sided
admissible function in $\Hnty\big(\Omin,\CB(Y,U)\big)$. Suppose that the
following conditions hold:

\nr{1} $\CO_{\A,\C}:X\to E^2(\Omext,Y)$ and $\CO_{\A^*,\B^*}:X\to E^2(\Omext,U)$
are bounded injective operators;


\nr{2} $\Phi$ corresponds to $\de$.

\noindent Then the observation
systems $(\A,\C)$ and
$(\A^*,\B^*)$ are dual with respect to $\de$.
\end{vartheorem}

\begin{rem} \enspace \enspace 
In \cite{YmodAandA},
we gave a wider definition of the
correspondence between $\Phi$ and $\de$.
It was required there that
$\tau$ and $\Phi$ belong to
wider functional classes than
classes $\Hnty$.
Theorem 9.5 from \cite{YmodAandA}
gives \emph{a necessary and sufficient condition} for
the duality of observation systems $(\A,\C)$ and
$(\A^*,\B^*)$, and in this sense the definition
in \cite{YmodAandA} is the adequate one.
For our purpose, the above formulation in Theorem A will suffice.
\end{rem}

The proof of Theorem \ref{thm dual} will consist in checking
conditions (1) and (2) of Theorem A for the triple
$(A,B,-\ka C)$.
The transfer function of this system can be simply expressed as
$$
\Phi(z)=\ka C(A-zI)^{-1}B.
$$
Indeed, for $z\in \rho(A)$,
$(A-zI)^{-1}B$ is bounded from $X$ to $X_{\frac{\psi(t)}{|t|+1}}$
(see Lemma \ref{lem on siA and H}, (ii)),  and
$C$ is bounded from $X_{\frac{\psi(t)}{|t|+1}}$ to $X$.

\section{THEOREM \ref{thm dual}: DETAILS OF PROOF}
\begin{lemma}
\label{lem Phi Hnty}
The transfer function $\Phi$ of the system
$(A,B,-\ka C)$ belongs to
\linebreak
$\Hnty(\Omext, \CB(X))$.
\end{lemma}

\begin{proof}                  
Since $A-zI$ is bounded from $X_{\psi^{-1}}$ to
$X_{\psi}$, $(A-zI)^{-1}$
is bounded from $X_{\psi}$ to $X_{\psi^{-1}}$ for
$z\in \rho(A)$.
It follows from \eqref{frm for res A} that
$$
\Phi(z)=\ka CB(\Ao-zI)^{-1}\cdot
B^{-1}[I+L(\Ao-zI)^{-1}]^{-1}B.
$$
By \eqref{norm Omext},
$B^{-1}[I+L(\Ao-zI)^{-1}]^{-1}B$ is an
$\Hnty(\Omext, \CB(X))$ function and
by \eqref{le 1r'},
$CB(\Ao-zI)^{-1}$ is an
$\Hnty(\Omext, \CB(X))$ function.
\end{proof}

We will need a domain
$$
\Omin' \defin \Om^\text{int}_{\mu-\si,R-\si}
$$
where
a small parameter $\si>0$ is chosen in such a way that
\eqref{n FC resolv B} and \eqref{norm L by resolv} still hold
true for $z\in \BC\sm \clos\Omin'$ for some $\eps>0$.
Then $\Phi$ belongs to $\Hnty(\BC\sm \clos\Omin', \CB(X))$.

Put $\phi_*:\BR\to\BR$ to be the even continuation of $\phi$ if
$\phi$ is defined on $[0,\infty)$, and
$\phi_*=\phi$, if $\phi$ is already defined on $\BR$.

\begin{lemma}
\nr{i} For all $s>1$ and $t\in \BR_+$, $\phi(st)\le s\phi(t)$;

\nr{ii} $\phi_*(s)\le (|s|+1)\phi_*(1)$ for $s\in \BR$;

\nr{iii} $\phi(s+t)\le \phi(s)+ \phi(t)$ for $s,t\ge0$.
\end{lemma}

\begin{proof}
(i) and (iii) follow from the concavity of $\phi$
on $[0,\infty)$ and the fact that $\phi(0)\ge0$.
Conditions on $\psi$ imply that $\psi$ grows on
$[0,\infty)$.
Therefore (ii) is obtained by putting
$t=1$ in (i).
\end{proof}

It is clear that there is $\mu_1>\mu$ such that
$$
\Om^{\text{int}}_{*,\,\mu_1}\supset \clos \Om^{\text{int}}_{\mu,R},
$$
where
$$
\Om^{\text{int}}_{*,\,\mu_1}\defin
\big\{z=x+iy\in\BC: \quad
x\in\BR,\,|y|<\mu_1\phi_*(x)
\big\}.
$$
We choose a real number $\ell>\|F\|$ and put
\beqn
\label{def ka0}
\ka_0=\ell +\mu_1\big(2+\al+\phi(\al)\big), \qquad \text{where} \quad \al=1+\ell \phi(1).
\neqn
Fix any $\ka>\ka_0$. Our aim is to prove that for any such $\ka$,
the conclusions of Theorems \ref{thm observ} and \ref{thm dual} hold.

\begin{lemma}
For any $x,t\in \BR$,
\beqn
\label{ineq x t phi}
\big|
t+i\ka \phi_*(t)-
\big(
x+i\mu_1 \phi_*(x)
\big)
\big|
\ge \ell \phi_*(t)
\neqn
\end{lemma}

\begin{proof}
Note that $\ka>\mu_1$.
Take any $x,t\in \BR$.
We distinguish two cases.

\nr{i} Suppose that  $|x|<\al(|t|+1)$. Then, by the previous Lemma,
$$
\phi_*(x)\le\al\phi_*(t)+\phi_*(\al),
$$
which gives that
\begin{multline*}
|\ka \phi_*(t)-\mu_1\phi_*(x)|
\ge
(\ka-\mu_1)\phi_*(t)-\mu_1
\big(\phi_*(t)+\phi_*(x)\big)   \\
\ge
(\ka-2\mu_1-\mu_1\al)\phi_*(t)-\mu_1\phi_*(\al)
\ge \ell \phi_*(t).
\end{multline*}
The last inequality is due to the facts that
$\ka>\ka_0$ and $\phi_*(t)\ge1$. Now
\eqref{ineq x t phi} follows.

\nr{ii}
Suppose that  $|x|\ge\al(|t|+1)$.
Since $\al-1=\ell\phi(1)$, we get
%
$$
|t-x|\ge|x|-|t|\ge(\al-1)|t|+\al>               
\ell\phi(1)|t|+\ell\phi(1)\ge \ell \phi_*(t),
$$
and this again gives \eqref{ineq x t phi}.
\end{proof}

\begin{lemma}
For any $z\in \Om^\text{int}_{*,\,\mu_1}$ and
any $t\in \cD(\phi)$,
\beqn
\label{ineq ka phi}
|t+i\ka \phi(t)-z|\ge \ell \phi(t).
\neqn
\end{lemma}

\begin{proof}
Since $\ka>\mu_1$, $t+i\ka \phi(t)$ is outside
$\clos\Om^\text{int}_{*,\,\mu_1}$. Hence the straight line
interval with endpoints in
$t+i\ka \phi(t)$ and $z$ contains a boundary point of
$\Om^\text{int}_{*,\,\mu_1}$, which has a form
$x+i\mu_1 \phi_*(x)$ for some $x$. Therefore
\eqref{ineq ka phi} follows from \eqref{ineq x t phi}.
\end{proof}



\

\begin{lemma}
\label{lem about de}
\nr{i}
$\Aka-zI+i\phi(\Ao)F$ is an isomorphism from
$X$ onto $X_{|t|+1}$ for all $z\in \Om^\text{int}_{\mu_1}$;

\nr{ii} $\de\in \Hnty \big(\Om^\text{int}_{\mu_1}, \CB(X)\big)$;

\nr{iii} For any $z\in \Om^\text{int}_{*,\,\mu_1}\sm\Omin'$,
$\de(z)$ is invertible, and
$$
\de^{-1}(z)=I+\Phi(z).
$$
For these values of $z$, the norms of $\de^{-1}(z)$ are uniformly bounded.
In particular, the norms $\|\de^{-1}(\cdot)\|$ are uniformly bounded on
$\pt\Omin$.
\end{lemma}

\begin{proof}
(i)
Since $\Aka-zI$ is an isomorphism from $X$ to $X_{|t|+1}$ for
$z\in \Om^\text{int}_{*,\,\mu_1}$, one has to check that
$I+i\phi(\Ao)(\Aka-zI)^{-1}F$ is an invertible operator on $X$.
By \eqref{ineq ka phi},
$$
\|
\phi(\Ao)(\Aka-zI)^{-1}F
\|\le \|F\| \sup_{t\in\cD(\phi)}
\frac{\phi(t)}{t+i\ka\phi(t)-z}
\le \frac {\|F\|} \ell<1
$$
for $z\in\Om^\text{int}_{*,\,\mu_1}$.
We obtain that for $z\in\Om^\text{int}_{*,\,\mu_1}$,
\beqn
\label{***}
\|
\big[I+i\phi(\Ao)(\Aka-zI)^{-1}F\big]^{-1}
\|
\le
\frac 1 {1-\|F\|/ \ell }<\infty,
\neqn
and our assertion follows.

(ii) It is easy to check that
\beqn
\label{*2}
\de(z) = I-i\ka
\big[I+i\phi(\Ao)(\Aka-zI)^{-1}F\big]^{-1}
(\Aka-zI)^{-1}\phi(\Ao).
\neqn
Since $\sup_{z\in \Om^\text{int}_{*,\,\mu_1}}\|(\Aka-zI)^{-1}\phi(\Ao)\|<\infty$,
the assertion follows from \eqref{***}.

(iii) By \eqref{def H z},
$$
\de(z)=\big[\Aka-zI+i\phi(\Ao)F\big]^{-1}(\Ao-zI)H(z).
$$
For $z\in \Om^\text{int}_{*,\,\mu_1}\sm\Omin'$,
$(\Ao-zI)H(z)$ is an isomorphism from $X$ to $X(\Ao)$
(see Lemma \ref{lem on siA and H}, (iv)), and
$\big[\Aka-zI+i\phi(\Ao)F\big]^{-1}$ is an isomorphism from
$X(\Ao)$ to $X$. Hence for these $z$, $\de(z)$ is invertible,
and
\beqn
\label{I+Phi}
\begin{aligned}
\de(z)^{-1}& =
\big[\Ao-zI+i\phi(\Ao)F\big]^{-1}
\big[\Aka-zI+i\phi(\Ao)F\big]            \\
& = I+ i\ka \big[\Ao-zI+i\phi(\Ao)F\big]^{-1} \phi(\Ao) \\
& = I+ i\ka \psi(\Ao)\big[A-zI\big]^{-1} \psi(\Ao) \\
& = I +\Phi(z).
\end{aligned}
\neqn
The definition of $\Omin'$ implies that $\Phi$ is an operator-valued
$\Hnty$ function in $\BC\sm\clos\Omin'$. Hence
$\|\de^{-1}(z)\|\le K<\infty$ in
$\Om^\text{int}_{\mu_1}\sm\Omin'$.
\end{proof}

\begin{proof}[Proofs of Theorems \ref{thm dual} and \ref{thm observ}]
Let $R$, $\ka_0$ be chosen as above, and take any
$\ka>\ka_0$. Put $\de=\de_\ka$.
We check the hypotheses of Theorem A for the system
$(A,B,-\ka C)$.
By Lemma \ref{lem on siA and H}, (v), $\CO_{A,C}$ is bounded.
Since $C=i\psi(\Ao)$ and $\psi\ne0$ on $\si(\Ao)$,
$\CO_{\Ao,C}$ is injective.
Statements (iii) and (iv) of Lemma \ref{lem on siA and H} imply that
$\CO_{A,C}$ is injective. By symmetry, the same can be said about
$\CO_{A^*,B^*}$.

By Lemma \ref{lem about de}, (ii) and (iii), $\de$ is two-sided admissible.
By \eqref{I+Phi},
$\Phi$ corresponds to $\de$, with $\tau(z)\equiv I$. Hence all the hypotheses
of Theorem A are fulfilled. Therefore observation systems
$(A,-\ka C)$ and $(A^*, B^*)$ are dual to each other
with respect to $\de$. In particular, $\CO_{A,C}$ is an isomorphism.

Statements (1) and (2) of Theorem \ref{thm observ} have been already verified.
Now statement (3) of this theorem follows
immediately from  Theorem 3.1 and Proposition 5.1 in \cite{YmodAandA}.
\end{proof}

It is easy to see that the same results hold for $\ka<-\ka_0$.

\section{THE CONTROL MODEL}
\label{sect contr mod}

In this Section, we give a dual formulation of our result in terms of
what we call the control model. Let us give first a control theory motivation.

Suppose we are given an abstract control system
$(\A,\B)$ such that
$\A$ is a generator of a bounded $C_0$ semigroup.
Associate with the pair $(\A,\B)$ the linear system
$$
\dot x(t)=\A x(t)+\B u(t).
$$
Consider the mapping $\cC_{\A,\,\B}$, which sends
an input $u\in L^2(\BR_-,U)$ to
$x(0)$. We define it first for smooth functions $u$ with compact support,
assuming that $x(t)=0$ for large negative times. For these functions,
$\cC_{\A,\,\B}$ is always well-defined.
Next we make an assumption that $\cC_{\A,\,\B}$ extends continuously to
$L^2(\BR_-,U)$. In the control theory, a system with
this property is called \textit{infinite time admissible}
\cite{Sta-libro}. Denote the extended map by the same symbol
$\cC_{\A,\,\B}$.
For an admissible system, define the controllability map
$$
\WABSTR :\EinU\to X
$$
by taking the composition map in the diagram
$$
\EinU
\underset{\cL^{-1}}{\longrightarrow} L^2(\BR_-,U)
\underset{\cC_{\A,\,\B}}{\longrightarrow} X
$$
where $\cL^{-1}$ is the inverse Laplace transform
and $\Omin=\Pi_-$.
By the usual convention, $E^2(\Pi_-)=H^2(\Pi_-)$.
We put  $\Omext=\Pi_+$. It is easy to see that
$$
\WABSTR(z-\la)^{-1}u=(\A-\la I)^{-1} \B u,
\qquad \la\in\Omext, \; u\in U.
$$

Let us return to the general case when
$(\A,\B)$ is an arbitrary abstract control system
and $\Omin$ an arbitrary admissible domain such that $\si(\A)\subset \clos\Omin$.
We take the last formula as a starting point of the
general definition of the transform
$W_{\A,\,\B}$.

Consider the linear set $\cH$ of rational holomorphic functions
from $\Omin$ to $U$ that are representable as finite sums
$$
f(z)=\sum_j  (z-\la_j)^{-1}\,u_j
$$
with $u_j\in U$ and $\la_j\in \Omext$. For each such $f\in \cH$, we put
$$
\WoABSTR f\,\defin\, \sum_j (\A-\la_jI)^{-1}\B u_j.
$$
It is easy to prove that $\cH$ is a dense subset of $\EinU$.
\begin{definitions}
\nr {1} Abstract control system $(\A,\B)$ is called \emph{admissible
with respect to} $\Omin$ if
$\WoABSTR$ extends to a continuous operator
$$
\WABSTR: \EinU \to X.
$$

\nr {2} Abstract control system $(\A,\B)$ is called \emph{exact
with respect to} $\Omin$ if
it is admissible with respect to this domain and the image of the extended map
$\WABSTR$ is the whole space $X$.
\end{definitions}

It is easy to see that $\WABSTR$ splits the multiplication operator by $z$
on %
%
%
\linebreak
$\EinU$ with the operator $\A$; more exactly,
\begin{equation}
\label{qA}
\WABSTR [q(z)f(z)]=q(\A)\big(\WABSTR f\big), \qquad f\in \EinU
\end{equation}
for any rational scalar function \enspace $q\in H^\infty(\Omin)$.
For these functions, $q(\A)$ is bounded. It follows from this equation that
$\ker\WABSTR$ is invariant under the multiplication by
rational functions in $\Hnty(\Omin)$.

For any admissible abstract control system
$(\A,\B)$, we define
the quotient operator $\WhABSTR$ by factoring
$\WABSTR$ by its kernel:
\begin{align*}
&\WhABSTR:\EinU\big/\ker\WABSTR\to X, \\
&\WhABSTR \big(f+ \ker\WABSTR) \defin \WABSTR f.
\end{align*}
By the Beurling--Lax--Halmos theorem, $\ker\WABSTR$ has a form
\beqn
\label{ker eqs de}
\ker\WABSTR=\de \EinY
\neqn
for a Hilbert space $Y$ and an admissible function
$\de\in \Hnty\big(\Omin,\CB(Y,U)\big)$.
\begin{definition}
Any admissible function $\de$
satisfying \eqref{ker eqs de}
will be  called
\textit{a generalized characteristic function } of abstract control
system $(\A,\B)$.
\end{definition}

Put
(as before)
$\bomin=\{\bar z: z\in\Omin\}$ and $\bomext=\{\bar z: z\in\Omext\}$.
Notice that the pairing
\beqn
\label{cauchy pair}
\langle
f,g
\rangle\defin
\frac 1 {2\pi i}
\int_{\Ga}
\langle
f(z),g(\bar z)
\rangle\, dz,
\quad f\in \EinU,\, g\in \BEoutU
\neqn
defines a duality between Hilbert spaces
$\EinU$ and $\BEoutU$, see \cite{YmodAandA}.

\begin{lemma}
\label{rel ctr obs mods}
\nr{1} Abstract control system $(\A,\B)$ is admissible
with respect to $\Omin$
if and only if abstract observation system $(\A^*,\B^*)$ is admissible
with respect to $\bomin$;

\nr{2} If any of the above two assertions holds, then
$$
\CO_{\A^*,\,\B^*}=W_{\A,\,\B}^{*}
$$
with respect to the pairing \eqref{cauchy pair}.
\end{lemma}

\begin{proof}
The Cauchy pairing \eqref{cauchy pair} extends to rational functions
$f\in \cH$ and \textit{all} holomorphic functions
$g:\Omext\to U$ by putting
$$
\langle (z-\la)^{-1}u, g\rangle\defin \langle u, g(\la)\rangle
$$
and extending this formula by linearity. For functions $g\in\BEoutU$, it
is the same pairing. It is plain to check the identity
$
\langle f, \CO_{\A^*,\,\B^*}x\rangle=\langle \WoABSTR f, x\rangle
$
for $f\in\cH $ and $x\in X$. Both statements of Lemma are
easy consequences of this formula.
\end{proof}

\begin{lemma}
\label{lem equiv ctr obs}
Suppose we are given an abstract control system $(\A, \B)$, which is admissible
with respect to an admissible domain $\Omin$ and
a two-sided admissible function
$\de\in \Hnty\big(\Omin,\CB(Y,U)\big)$.

Then the following are equivalent.

\nr{1} $\CO_{\A^*,\,\B^*}:X\to \hdest$
is an isomorphism;

\nr{2} $(\A,\B)$ is an exact abstract control system and
$\ker\WABSTR=\de \EinY$;

\nr{3} $\ker\WABSTR=\de \EinY$ and $\WhABSTR$ is an isomorphism.
\end{lemma}

\begin{proof}
If follows from Lemma \ref{rel ctr obs mods} that
the closure of the range of $\CO_{\A^*,\,\B^*}$
in $\EoutU$ equals to the annihilator of $\ker \WABSTR$ with respect to
the Cauchy pairing \eqref{cauchy pair}. By the Banach theorem, the range
of $\CO_{\A^*,\,\B^*}$ is closed if and only if
$\WABSTR=\CO_{\A^*,\,\B^*}^*$ is onto. Finally,
$$
\hdest=\big(\de\EinY\big)^\perp
$$
(see \cite[Proposition 2.5]{YmodAandA}). These remarks imply the
equivalence of statements (1)--(3).
\end{proof}

\begin{definition}
For a domain $\Omin$, Hilbert spaces
$U$, $Y$ and a fixed two-sided admissible function
$\de\in \Hnty(\Omin,\CB(Y,U))$,
we consider \emph{the control model space}, which is the quotient space
\beqn
\EinU\big/\de \EinY.
\neqn
For a function $f$ in $\EinY$, we put $[f]=f+\de \EinY$ to be its
coset in this quotient space.

\textit{The model operator $\whMz$} on this space is simply the quotient
operator of multiplication by the independent variable $z$. It is given by
\begin{align*}
\cD(\whMz)&=
\big\{
[f]: \; f, M_zf\in \EinY
\big\},                                \\
\whMz [f] &\defin [zf], \qquad \text{ if }
[f]\in \cD(\whMz)\; \text{ and } f, zf\in \EinY.
\end{align*}
\end{definition}

\begin{theorem}
\label{last thm}
Let $A$ be an operator given by \eqref{intro def A},
\eqref{def A}, where $\psi$ is a function satisfying
(1)--(6). Put $B=\psi(A)$, and define
$\de_\ka$  as in Theorem \ref{thm observ}.
Then there exist  $R>0$ and $\ka\in \BR$
such that for the corresponding
function $\de= \de_\ka$ and for the domains
$\Omin$, $\Omext$, given by
\eqref{Omin Omext}, the following are true.

\nr{1}  $(A,B)$ is an exact control system in $\Omin$ and $\de$
is its generalized characteristic function (that is,
$\ker W_{A,B} =\de \EinX$);

\nr{2} Operator
$$
\wh W_{A,B}: \EinX\big/\de \EinX \to X
$$
is an isomorphism that transforms operator $A$ into the
quotient multiplication operator
$\whMz$ on the control model space.
In particular,
$\wh W_{A,B}\,\cD(\whMz)=\cD(A)$, and
\beqn
A\,\wh W_{A,B} [f]= \wh W_{A,B}\,\whMz [f], \qquad
\forall f\in \cD(\whMz).
\neqn
\end{theorem}

\begin{proof}
Take $\mu$, $R$ and $\ka_0$ as in Sections 1--4, and put
$\Omin=\Om^\text{int}_{\mu,R}$. Take any $\ka>\ka_0$, and
put $\de=\de_\ka$. By
Theorem \ref{thm dual}, observation system $(A^*, B^*)$ is exact
and $\CO_{A^*,\, B^*}$ is an isomorphism from $X$ onto
$\hdest$. Hence, by Lemma \ref{lem equiv ctr obs},
(1) holds, and
$\wh W_{A,B}$ is an isomorphism.
The splitting properties of this isomorphism, stated
in (2),  follow easily from
\eqref{qA}.
\end{proof}

The calculus $q\mapsto q(A)\in \CB(X)$ is defined for any
rational function $q$ in $\Hnty(\Omin)$.

\begin{corollary}
\label{cor}
The above functional calculus $q\mapsto q(A)$
extends by continuity to
an $\Hnty(\Omin)$ functional calculus for $A$. In particular,
\beqn
\label{K-spectr}
\|f(A)\|\le K\sup_{z\in\Omin}|f(z)|
\neqn
for
$f\in\Hnty(\Omin)$, where
$K=\|\wh W_{A,B}\|\,\|\wh W_{A,B}^{-1}\|$.
\qed
\end{corollary}

\begin{corollary}
\label{cor2}
Operator $A$ admits \emph{a skew normal dilation} on $\pt \Omin$ in the
following sense. There exists a Hilbert space $K$ and an unbounded
operator $N$, acting on $K$, that has the following properties.

\nr{i} $N$ is similar to an unbounded normal operator;

\nr{ii} $\si(N)$ is contained
in $\pt\Omin$ and is absolutely continuous with respect to the arc length measure;

\nr{iii} $q(A)=P_X q(A)|X$ for any rational function $q$ in $\Hnty(\Omin)$.
\end{corollary}

\begin{remarks}
\nr{1} The quotient operator $\wh M_z$
of multiplication by $z$ can also be defined for the case of two-sided
admissible function $\de$ on a \emph{bounded}
domain $\Omin$; in this case $\wh M_z$ is bounded.
For the case of a simply connected $\Omin$, one can substitute
$\de$ by its inner part $\de_i$, which comes from a canonical factorization
$\de=\de_i\de_e$ in this model.
In particular, in the case
when $\Omin$ is the unit disc $\BD$, the model operator of Theorem
\ref{last thm} becomes exactly the Nagy--Foia\c{s} model.
In the general case of bounded or unbounded (simply connected)
admissible domain $\Omin$,
one can identify $\wh M_z$ with $\ga(T)$, where
$\ga:\BD\to \Omin$ is a conformal mapping and
$T$ is the Nagy--Foia\c{s} model operator, whose
Nagy--Foia\c{s} characteristic function is $\de_i\circ\ga$
(more precisely, the pure part of this function). We understand
$\ga(T)$ in the sense of the Nagy--Foia\c{s} theory.
See \cite[\S5]{YmodAandA} for more details.

\nr{2} Inequality \eqref{K-spectr} implies that
$\clos\Omin$ is a so-called $K$-spectral set of $A$.
As Pisier proved in 1997 (see \cite{Pisier}),
the fact that a set $T$ is $K$-spectral of an operator does not imply the existence
of a skew dilation of this operator
to a normal operator whose spectrum is contained in
$\prt T$.
We refer to
\cite{Sarason, Hildebrandt, Lange, Lewis, Paulsen}
 and others for more information on
$K$-spectral sets of operators and positive and negative results on
similarity.

The existence of invariant subspaces for
Banach space operators such that the unit disc is their $K$-spectral
set (with a certain additional condition)
has been proved in \cite{AmbrMull}.
In our situation, by applying the results
by Nagy and Foia\c{s}
 \cite[Chapter VII]{SzNF}, we can
describe all rationally
invariant subspaces of operator $A$
under consideration
in terms of regular factorizations
of $\de$. A rationally invariant subspace of $A$ is, by definition,
an invariant subspace of $(A-\la I)^{-1}$ for  all $\la\in\rho(A)$.

The results by Nagy and Foia\c{s} are also applicable
to the lifting of the commutant of $A$.

It would be interesting to use the results of
\cite{Kupin} to give a necessary and sufficient condition for
similarity of $A$ to a normal operator. We refer to
\cite{BenamNik} and \cite{KupinTreil} for additional information.

\nr{3} Suppose we have an unbounded operator $A=\Ao+L$, where
$L$ has been represented as $L=i\psi(\Ao)F\psi(\Ao)$, so that
conditions (1)--(6) on $\psi$ and $\Ao$ are fulfilled. Take
any function $\psi_1$ that satisfies (1)--(4) and such that
 $\psi_1\ge \psi$ on $\cD(\psi)$. Then
$L=i\psi_1(\Ao)F_1\psi_1(\Ao)$, where
$\|F_1\|_\text{ess}\le\|F\|_\text{ess}$, so that
conditions (1)--(6) are also fulfilled for $\psi_1$ and $\Ao$.
Hence, whenever our construction yields a model in some parabolic domain,
it also gives a model in larger parabolic domains, with
other auxiliary operators $B$ and $C$.
\end{remarks}

\section{PROOFS OF AUXILIARY LEMMAS}

\begin{proof}[Proof of Lemma \ref{lem un discs}]
Remind that $r'$ and $k$ were chosen so as to satisfy
\eqref{ineq r' k mu}.
There exists $t_0>0$ such that
\beqn
\label{**}
{\phi(t)}/t<k\qquad \text{for} \quad t\ge t_0.
\neqn
Let us prove that for these $t$, the disc
$B(t, r'\phi(t))$ is contained in
$$
\Om_+^\text{int}\defin \Om^\text{int}_\mu\cap \{z:\Re z>0\}.
$$
The vertical line $\Re z=t$ divides the disc
$B(t, r'\phi(t))$ into two halves.
First notice that the right half-disc is contained in
$\Om_+^\text{int}$. Indeed, if
$z=x+iy\in B(t, r'\phi(t))$ and $x\ge t$, then
$|y|\le r'\phi(t)\le \mu\phi(x)$. It remains to prove that the
left half-disc is contained in $\Om_+^\text{int}$.
Consider the triangle $T$ with vertices at
points $0$ and $\tau_\pm=t\pm i\mu\phi(t)$. Since
$\Om_+^\text{int}$ is convex, $\inter T\subset \Om_+^\text{int}$.
It is easy to check, using
\eqref{def mu0}, \eqref{**} and the second inequality in
\eqref{ineq r' k mu} that the distances from the point $t$ to the
sides $[0,\tau_\pm]$ of the triangle $T$ equal to
$$
{t\mu \phi(t)}\big/ {\sqrt{t^2+\mu^2\phi^2(t)}},
$$
which is greater than $r'\phi(t)$. This proves that the left
half of the disc is contained in $\Om_+^\text{int}$.
We conclude that
\eqref{B subset Om} holds for $|t|\ge t_0$.
The union of discs
$B(t, r'\phi(t))$, $|t|< t_0$ is bounded, and the
statement of Lemma follows.
\end{proof}

\begin{proof}[Proof of Lemma \ref{lem est int}]
Fix some $t_0>0$ such that
\eqref{**} holds. Put $\rho=1/(2k)$.
For a fixed $x$, divide $\cD(\phi)\sm (-2R,2R)$ into a countable union of sets
$I_n=I_n(x)$, $n\ge0$, where
\begin{align}
I_0&=\big\{
t\in \cD(\phi): |t|\ge 2R, \;|t-x|\le \rho\phi(x)
\big\}, \nn \\
I_n&=\big\{t\in \cD(\phi): |t|\ge 2R, \;
2^{n-1}\rho\phi(x)\le|t-x|\le 2^n\rho\phi(x)
\big\}, \quad n\ge 1.
\nn
\end{align}
Then $|I_n|\le 2^{n}\rho\phi(x)$ for $n\ge 1$.

Put
$$
\Ga_n'=\big\{
x+iy\in \Ga: x\in I_n
\big\}, \qquad n \ge 0,
$$
and
$\Ga''=
\big\{
x+iy\in \Ga: |x|\le 2R
\big\}$.
Then $\Ga=\big(\bigcup_{n\ge0} \Ga_n'\big)\cup \Ga''$.

We parametrize $\Ga'_n$ by $z(t)=t\pm i\mu\phi(t)$, $t\in I_n$.
We have $|dz(t)|/dt\le C_1$ on all curves $\Ga'_n$. For $n\ge 1$,
\begin{align}
\int_{\Ga_n'}\frac{|dz|}{|x-z|^2}\le 2C_1
\int_{I_n}\frac{dt}{|x-z(t)|^2} 
\le 2C_1
\int_{I_n}\frac{dt}{|x-t|^2}
\le \frac {2C_1|I_n|}{2^{2n-2}\rho^2 \phi(x)^2}\le
\frac {2^{3-n}C_1}{\rho \phi(x)}.
\label{n ge 1}
\end{align}
Next, for all $x>0$ sufficiently large,
$\phi$ increases on the interval
$[\frac x 2, 2x]$, which contains $I_0(x)$
(we use \eqref{**}).  Hence for
$t\in I_0(x)$,
$
\phi(t)\ge \phi(x/2)\ge \phi(x)/2.
$
Similar estimates hold for $x<0$ with large $|x|$, and we obtain that
\beqn
\label{n eq 0}
\int_{\Ga_0'}\frac{|dz|}{|x-z|^2}\le 2C_1
\int_{I_0}\frac{dt}{|x-z(t)|^2}
\le
\frac {8 C_1|I_0|}
{\mu^2\rho^2\phi(x)^2}
\le
\frac {16C_1}{\mu^2\rho\phi(x)}.
\neqn
By \eqref{n ge 1} and  \eqref{n eq 0},
$$
\int_{\Ga'}\frac{|dz|}{|x-z|^2}\le
\frac {16C_1}{\mu^2\rho\phi(x)}
+
\sum_1^{\infty}
\frac {2^{3-n}C_1}{\rho\phi(x)}
\defin\frac {C_2}{\phi(x)}.
$$
Since $\int_{\Ga''}
\frac{|dz|}{|x-z|^2}\sim
\frac {C_3}{|x|^2}\le \frac {1}{\phi(x)}
$
for large $|x|$, we obtain the statement of Lemma.
\end{proof}


\begin{thebibliography}{99}







\bibitem{Agranovich}
\textsc{M. S. Agranovich}, Elliptic operators on closed manifolds,
in Partial Differential Equations VI (Encyclopaedia Math. Sci.,
vol. 63), Springer-Verlag, Berlin---New York 1994, pp. 1--130.

\bibitem{AmbrMull}
\textsc{C. Ambrozie, V.  M\"uller}, Invariant subspaces for
polynomially bounded operators, \textit{J. Funct. Anal.} 213
(2004), no. 2, 321--345.


\bibitem{Arlinsky}
\textsc{Yu. Arlinsky}, Characteristic functions of maximal
sectorial operators, Recent advances in operator theory
(Groningen, 1998), 89--108, Oper. Theory Adv. Appl., 124,
Birkh\"user, Basel, 2001.

\bibitem{Bad_Crou_Del}
\textsc{C. Badea, M. Crouzeix, B.  Delyon}, Convex domains and
$K$-spectral sets. \textit{Math. Z.} 252 (2006), no. 2, 345--365.

\bibitem{BenamNik}
\textsc{N.-E. Benamara, N.  Nikolski}, Resolvent tests for
similarity to a normal operator. \textsc{Proc. London Math. Soc.}
(3) 78 (1999), no. 3, 585--626.



\bibitem{BirmSol}
\textsc{M. Birman, M. Solomjak}, \textsc{Spectral theory of
selfadjoint operators in Hilbert space},
Mathematics and its Applications (Soviet Series). D. Reidel Publ.
Co., Dordrecht, 1987.

\bibitem{Crouz_parb}
\textsc{M. Crouzeix}, Operators with numerical range in a
parabola. \textit{Arch. Math. (Basel)} 82 (2004), no. 6, 517--527.

\bibitem{Crouz_Del_sector}
\textsc{M. Crouzeix, B. Delyon}, Some estimates for analytic
functions of strip or sectorial operators. \textit{Arch. Math.
(Basel)} (2003), no. 5, 559--566.


\bibitem{CurZw}
\textsc{R. F. Curtain, H. J. Zwart}, \textit{An introduction to
infinite-dimensional linear systems theory}, Springer, New York,
1995.


\bibitem{DunfSchw2}
\textsc{N. Dunford,  J. T. Schwartz}, \textit{Linear operators.
Part II: Spectral theory. Self adjoint operators in Hilbert
space}. Interscience Publishers John Wiley \& Sons N.Y. -- London
1963.


\bibitem{Dur}
\textsc{P. Duren}, \textit{Theory of $H^p$-spaces}, Acad. Press,
N. Y., 1970.

\bibitem{Hildebrandt}
\textsc{S. Hildebrandt}, The closure of the numerical range of an
operator as spectral set, \textit{Comm. Pure Appl. Math.} 17
(1964),  415--421.

\bibitem{Keldysh}
\textsc{M. Keldy\v s}, \textit{The completeness of eigenfunctions
of certain classes of nonselfadjoint linear operators} (Russian),
Uspehi Mat. Nauk 26 (1971), no. 4(160), 15--41. Engl. transl.:
Russian Math. Surveys 26 (1971), no. 4, 15--44.

\bibitem{Kupin}
\textsc{S. Kupin}, Linear resolvent growth test for similarity of
a weak contraction to a normal operator. \textit{Ark. Mat.} 39
(2001), no. 1, 95--119.

\bibitem{KupinTreil}
\textsc{S. Kupin, S. Treil}, Linear resolvent growth of a weak
contraction does not imply its similarity to a normal operator.
\textit{Illinois J. Math.} 45 (2001), no. 1, 229--242.

\bibitem{Lange}
\textsc{R. Lange}, Essentially subnormal operators and
$K$-spectral sets. \textit{Proc. Amer. Math. Soc.} 88 (1983), no.
3, 449--453.

\bibitem{leMe3}
\textsc{C. Le Merdy}, Similarities of $\omega$-accretive
operators. International Conference on Harmonic Analysis and
Related Topics (Sydney, 2002), 84--95, \textit{Proc. Centre Math.
Appl. Austral. Nat. Univ.}, 41, Austral. Nat. Univ., Canberra,
2003.



\bibitem{Lewis}
\textsc{K. A. Lewis}, Intersections of $K$-spectral sets.
\textit{J. Oper. Theory} 24 (1990), no. 1, 129--135.



\bibitem{Markus}
\textsc{A. Markus}, \textit{Introduction to the spectral theory of
polynomial operator pencils},
With an appendix by M. V. Keldysh. Translations of Mathematical
Monographs, 71. American Mathematical Society, Providence, R.I.,
1988.

\bibitem{Mark_Mats}
\textsc{A. S. Markus, V. Matsaev}, Operators generated by
sesquilinear forms and their spectral asymptotics. Linear
operators and integral equations. \textit{Mat. Issled.} No. 61,
(1981), 86--106, 157. (in Russian)



\bibitem{Naboko}
\textsc{S. Naboko}, Functional model of perturbation theory and
its applications to scattering theory. (Russian) \textit{Boundary
value problems of mathematical physics}, 10. \textit{Trudy Mat.
Inst. Steklov.} 147 (1980), 86--114, 203.







\bibitem{NikBook2}
\textsc{N. Nikolski}, \textit{Operators, functions, and systems:
an easy reading.} Vol. 2. Model operators and systems.
Mathematical Surveys and Monographs, 93. American Mathematical
Society, Providence, R.I., 2002.

\bibitem{Paulsen}
\textsc{V. I. Paulsen}, \textit{Completely bounded maps and
dilations.} Pitman Research Notes in Mathematics Series, 146.
Longman Scientific $\&$ Technical, Harlow; John Wiley $\&$ Sons,
Inc., New York, 1986.


\bibitem{Pisier}
\textsc{G. Pisier}, \textit{Similarity problems and completely
bounded maps}. Second, expanded edition. Includes the solution to
"The Halmos problem". Lecture Notes in Mathematics, 1618.
Springer-Verlag, Berlin, 2001. viii+198 pp.



\bibitem{Put_Sund}
\textsc{M. Putinar, S. Sundberg}, A skew normal dilation of the
numerical range of an operator, \textit{Math.  Ann.} 331 (2005),
345--357.

\bibitem{Sarason}
\textsc{D. Sarason}, On spectral sets having connected component,
\textit{Acta Sci. Math (Szeged)} 26 (1965), 289--299.




\bibitem{Sta-libro}
\textsc{O. Staffans}, \textit{Well-posed linear systems.}
Encyclopedia of Mathematics and its Applications, 103. Cambridge
University Press, Cambridge, 2005.


\bibitem{SzNF}
\textsc{B. Sz.-Nagy,  C. Foias}, \textit{ Harmonic analysis of
operators on a Hilbert space}, North-Holland, Amsterdam, and Akad.
Kiad\'o, Budapest, 1970.


\bibitem{Tikh1}
\textsc{A. Tikhonov}, A functional model and duality of spectral
components for operators with a continuous spectrum on a curve.
\textit{Algebra i Analiz} 14 (2002), no. 4, 158--195; transl. in
\textit{St. Petersburg Math. J.} 14 (2003), no. 4, 655--682


\bibitem{Tikh2}
\textsc{A. Tikhonov}, Transfer functions for ``curved''
conservative systems. Recent advances in operator theory, operator
algebras, and their applications, 255--264, Oper. Theory Adv.
Appl., 153, Birkh\"user, Basel, 2005.

\bibitem{Tikh3}
\textsc{A. Tikhonov}, Factorizations and invariant subspaces for
weighted Schur classes, Arxiv math.FA/0510135

\bibitem{Triebel}
\textsc{H. Triebel}, \textit{Interpolation theory, function
spaces, differential operators.} Second edition. Johann Ambrosius
Barth, Heidelberg, 1995.

\bibitem{YmodAandA}
\textsc{D.V. Yakubovich}, Linearly similar model of Sz.-Nagy --
Foias type in a domain (in Russian), \textit{Algebra i Analiz},
15, No. 2 (2003), 180-227, English transl. in \textit{St.
Petersburg Math. J.} 15, No. 2 (2004), 289-321.





 \end{thebibliography}
\end{document}